# Eine nützliche Lösung für das Sammelbilderproblem

Niklas Braband[1], Sonja Braband[1] und Malte Braband[2]

[1]Gymnasium Neue Oberschule, Braunschweig, Deutschland
{niklas.braband|sonja.braband}@no-bs.de
[2]Technische Universität Braunschweig, Deutschland
m.braband@tu-braunschweig.de

**Zusammenfassung.** Das Sammelbilderproblem ist eines der wenigen mathematischen Probleme, die regelmäßig in den Schlagzeilen der Nachrichten vorkommen. Dies liegt einerseits an der großen Popularität von Fußball-Sammelbildern (manchmal als Paninimania bezeichnet) und andererseits daran, dass es bisher keine Lösung gibt, die alle relevanten Effekte wie Nachkaufen oder Tauschen berücksichtigen. Wir haben bereits in früheren Arbeiten nachgewiesen, dass die klassischen Annahmen nicht der Realität entsprechen. Deshalb stellen wir zunächst neue Annahmen auf, die die Praxis besser abbilden. Darauf aufbauend können wir Formeln für die mittlere Anzahl benötigter Bilder (sowie deren Standardabweichung) ableiten, die alle in der Praxis auftauchenden Effekte berücksichtigen. Damit können Sammler die mittleren Kosten eines Albums sowie deren Standardabweichung nur mit Hilfe von elementaren Rechnungen bestimmen. Für praktische Zwecke ist das Sammelbilderproblem damit gelöst.

**Abstract.** The Coupon Collector's Problem is one of the few mathematical problems that make news headlines regularly. The reasons for this are on one hand the immense popularity of soccer albums (also known as Paninimania) and on the other hand that so far no solution is known that is able to take into account all effects such as replacement (limited purchasing of missing stickers) or swapping. In previous papers we have proven that the classical assumptions are not fulfilled in practice. Therefore we define new assumptions that match reality better. Based on these assumptions we are able to derive formulae for the mean number of stickers needed (and the associated standard deviation) that are able to take into account all effects that occur in practical collecting. Thus collectors can estimate the average cost of completion of an album and its standard deviation just based on elementary calculations. From a practical point of view we consider the Coupon Collector's problem as solved.

# 1    Einführung

Sammelbilder gibt es zu vielen Themen, z. B. zuletzt Fußballbilder zur EM 2016. Und es gibt viele Diskussionen über sie, z. B. zu viele Doppelte, manche Bilder scheinen seltener vorzukommen, einige Sammler vermuten sogar Betrug der Hersteller. Das Sammelbilderproblem ist eines der wenigen mathematischen Probleme, die regelmäßig in den Schlagzeilen der Nachrichten vorkommen. Dies liegt einerseits an der großen Popularität von Fußball-Sammelbildern (manchmal als Paninimania bezeichnet) und andererseits daran, dass es bisher keine Lösung gibt, die alle relevanten Effekte wie Nachkaufen oder Tauschen berücksichtigen.

Wikipedia [1] schreibt "Das Sammelbilderproblem .... befasst sich mit der Frage, wie viele zufällig ausgewählte Bilder einer Sammelbildserie zu kaufen sind, um eine komplette Bildserie zu erhalten". D. h., es geht darum, eine bestimmte Anzahl B von Bildern zu sammeln, um ein Album zu vervollständigen. Man kauft die Bilder nicht einzeln, sondern in Päckchen zu je P Bildern. Der Hersteller bietet in der Regel an, dass jeder Sammler einmalig eine begrenzte Anzahl von Bildern K (zu einem erhöhten Preis) nachkaufen kann, um seine Sammlung zu vervollständigen. Bei Panini war z. B. in Deutschland zur EM 2016 B=680, P=5 und K=50. Weiter ist noch der Preis p eines Päckchens wichtig, sowie der Preis b für die Nachbestellung eines Bildes (im Beispiel p=70 Cent; b=20 Cent, ohne Berücksichtigung von Nebenkosten wie Porto).

Weiter ist es noch wichtig, ob man alleine sammelt oder Bilder tauscht, z. B. in Tauschbörsen oder mit Freunden. F soll die Anzahl der Freunde sein (einschließlich des Sammlers selbst).

Klassisch werden im Sammelbilderproblem die folgenden Annahmen getroffen:

> K1. Die Bilder sind zufällig auf die Päckchen verteilt. D. h. der Hersteller mischt bei der Herstellung gut.
> K2. Alle Bilder kommen gleich häufig vor. D. h. der Hersteller betrügt nicht durch absichtliche Verknappung von Bildern.
> K3. Es wird in einer Sammelgemeinschaft fair getauscht, d. h. ein Bild gegen ein anderes.
> K4. Es gibt keine Rabatte, alle Bilder (außer Nachbestellungen) sind gleich teuer.
> K5. In einem Päckchen kommt kein Bild doppelt vor.

Aber man erkennt schnell, dass Sammeln ein teures Vergnügen ist. Denn um das Panini EM-Album zu füllen, hätte man mindestens 136 Päckchen zu 70 Cent kaufen müssen, was 95,20€ gekostet hätte.

Für den Handel werden Großpackungen, auch Display genannt, angeboten z. B. mit D=500 Bildern in 100 Päckchen. Das ist meistens wesentlich billiger, als die Päckchen einzeln zu kaufen. Dabei geht man häufig davon aus, dass in einer Box alle Bilder verschieden sind, oder zumindest wesentlich weniger Doppelte vorkommen. Wir hatten in einer früheren Arbeit [2] nachgewiesen, dass aufgrund von systematischen Abweichungen im Herstellungsprozess der Päckchen wirklich weniger Doppelte vorkommen als bei zufälliger Mischung. Aber es gibt nur wenige Serien, bei denen es keine Doppelten im Display gibt und die Fußballbilder zählen nicht dazu.

Es besteht weitgehende Einigkeit darüber [1], dass die folgende Strategie optimal ist, um das Album zu füllen:

1. Kaufe eine Box mit D Bildern.
2. Kaufe zusätzliche Päckchen und tausche so viele Doppelte wie möglich, bis maximal K Sticker in der Sammlung fehlen.
3. Kaufe die fehlenden K Bilder beim Hersteller nach.

Allerdings konnten die Kosten für diese Strategie nur durch Simulation [3] bzw. durch Näherungsformeln [4] bestimmt werden. Dies liegt vor allem daran, dass in der klassischen Annahme K3 eine sehr spezielle Annahme über das Tauschen gemacht wird [5], nämlich dass die Sammelgemeinschaft so lange weitersammelt, bis alle Alben der Sammler gefüllt sind. Dies entspricht dem Fall, dass ein Sammler F Alben sammeln würde und führt zu komplexen Abhängigkeiten.

In der Literatur sind viele Ergebnisse über das klassische Sammelbilderproblem [1, 3, 5, 6] zu finden. Nach [1] gilt für die mittlere Anzahl M zu kaufender Bilder ohne Berücksichtigung von K5 bei einem Sammler $BH(B)$, wobei $H(B)$ die harmonische Summe mit B Summanden bedeutet. D. h. der Faktor f =M/B beträgt hier gerade $H(B)$. Außerdem ist dort angegeben, dass die harmonische Summe gut durch den natürlichen Logarithmus angenähert werden kann, d. h. für große B gilt näherungsweise f=ln(B). Für das EM-Album würde dies bedeuten, dass ein Einzelsammler etwa 678€ ausgeben müsste.

Dieselbe Herleitung unter Berücksichtigung von K5 ist wesentlich komplizierter [6] und führt für das EM-Album nur zu einer geringen Ersparnis von im Mittel ca. 12 Bildern, sodass in der Regel K5 nicht berücksichtigt wird. Im engeren Sinn steht K5 auch im Widerspruch zu K1.

Vor der letzten EM verbreitete sich noch eine alternative Herleitung [7] rasant im Internet. Dabei wurde behauptet, dass unter K5 der Faktor nur $H(B/P)$ betragen würde und damit das EM-Album nur etwa 523€ kosten würde. Aber dieses Ergebnis erwies sich als falsch [8]. Außerdem lieferte es wie viele klassische Ergebnisse viel zu hohe Kosten, da das Nachkaufen und Tauschen nicht berücksichtigt wurden.

Für das Tauschen nach K3 ist bekannt [5], dass die mittlere Anzahl zu kaufender Bilder bei F Sammlern näherungsweise

$$B \ln(B) + B(F - 1)\ln(\ln(B)) \quad (1)$$

beträgt. Dies bedeutet, dass für eine große Anzahl F von Sammlern der Faktor f= M/(BF) gegen ln(ln(B)) strebt. Für das EM-Album wären das etwa 178€.

Der große Nachteil an diesen sowie allen anderen bisher veröffentlichten Ansätzen ist, dass sie die Realität nicht treffen. Einerseits stimmen die Annahmen nicht mit der Realität überein und andererseits bilden die Ansätze nicht alle Parameter angemessen ab, sodass überzogene Kosten resultieren. Deswegen werden wir in einem ersten Schritt praxisgerechte Annahmen ableiten.

## 2    Diskussion der Annahmen

Wir diskutieren jetzt die Gültigkeit der klassischen Annahmen, um unsere Sammelalbumformel auf möglichst realistischen Annahmen aufzubauen. Auch in einem DMV-Beitrag [9] wurde schon begonnen, u. a. auf Grundlage unserer Arbeiten, die klassischen Annahmen, insbesondere K1 und K2, in Frage zu stellen.

Die Annahme K1, dass die Bilder zufällig auf die Päckchen verteilt sind, d. h., der Hersteller mischt bei der Herstellung ordentlich, konnten wir in unserer letzten Arbeit widerlegen [2]. Es gab Abweichungen beim Herstellungsprozess. Die Verpackungsmaschine Fifimatic mischt die Bilder nicht zufällig. Allerdings fällt dies dem Käufer am Kiosk, der überwiegend einzelne Päckchen kauft, nicht auf, sondern nur dem Käufer eines Displays. Er erhält im Mittel deutlich weniger Doppelte als bei zufälliger Mischung. Das konnten wir sowohl bei unseren eigenen Untersuchungen beobachten, als auch bei Rezensionen von Käufern bei Amazon. Dort wurden durchschnittlich 47 Doppelte berichtet [10], bei zufälliger Mischung hätte man 146 erwartet [2], was aber in keinem der 38 berichteten Fälle erreicht worden ist. Deswegen ändern wir die Annahme K1 wie folgt ab:

A1. Die Bilder sind nicht zufällig auf die Päckchen verteilt, daher kommen in einem Display im Mittel weniger Doppelte vor als bei zufälliger Mischung.

Der DMV-Beitrag stimmt hier zu [9], allerdings ist die Annahme etwas anders formuliert, nämlich dass alle Päckchenkombinationen gleich häufig vorkommen.

Die Annahme K2, dass alle Bilder gleich häufig vorkommen, d. h., der Hersteller betrügt nicht, z. B. durch absichtliche Verknappung von Bildern, hat sich bei unseren Untersuchungen bisher bestätigt. Nur bei Trading Cards kommen geplant Bilder seltener vor [1], aber das wäre ein anderes Thema. Es gab zur letzten Fußball-WM Berichte [11], dass manche Bilder in Deutschland seltener vorkommen, aber dies gilt nicht in allen untersuchten Ländern wie z. B. der Schweiz. Der DMV-Beitrag zitiert dieselbe Quelle. Allerdings könnten diese vereinzelten Abweichungen auf Produktionsfehler zurückzuführen sein, was nach unseren Untersuchungen des Produktionsprozesses plausibel wäre. Daher behalten wir diese Annahme bei.

Die Annahme K3, dass innerhalb einer festen Sammelgemeinschaft fair getauscht wird, d. h. ein Bild gegen ein anderes, scheint nur noch vereinzelt zuzutreffen. Sie ist allerdings bisher die Standardannahme bei den meisten Analysen [1, 5, 7]. Wir haben nur eine einzige Quelle gefunden, die davon abweicht [12]. Nach unseren Beobachtungen sowie den Medienberichten wird vor allem bei Tauschbörsen getauscht, entweder vor Ort, z. B. in der Schule oder im Kaufhaus oder Online über das Internet. Das heißt, diese Annahme entspricht nicht mehr der Realität und wir ersetzen sie durch:

A3: Es wird mit Hilfe von Tauschbörsen fair getauscht, d. h. ein Bild gegen ein anderes.

Auch die Annahme K4, dass es keine Rabatte gibt, d. h. alle Bilder (außer Nachbestellungen) sind gleich teuer, ist in der Praxis nicht erfüllt, denn Displays sind billiger als der Kauf von einzelnen Päckchen. Deswegen ersetzen wir sie durch:

A4: Die Bilder sind bei Nachbestellungen teurer und im Display billiger als beim Kauf einzelner Päckchen.

Die Annahme K5, dass in einem Päckchen kein Bild doppelt vorkommt, konnten wir in unserer letzten Arbeit nachweisen [2]. Sie steht allerdings im Widerspruch zu K1 und sorgt für die beobachteten Abweichungen, die zu A1 führen. In der Praxis bringt sie aber nur einen geringen Vorteil für den Sammler, z. B. maximal zwölf Bilder von 4828 Bildern beim EM-Album [1]. D. h. es ist eher eine gut klingende Werbeaussage. Allerdings werden die Berechnungen bei Berücksichtigung des Päckchen-Effekts sehr kompliziert, d. h. der Mehraufwand ist im Vergleich zur Verbesserung der Genauigkeit unverhältnismäßig groß, sodass wir diese Annahme streichen.

Zusammengefasst sehen unsere neuen Annahmen so aus:

A1: Die Bilder sind nicht zufällig auf die Päckchen verteilt, daher kommen in einem Display im Mittel weniger Doppelte vor als bei zufälliger Mischung.

A2: Alle Bilder kommen gleich häufig vor, d. h., der Hersteller betrügt nicht, z. B. durch absichtliche Verknappung von Bildern.

A3: Es wird mit Hilfe von Tauschbörsen fair getauscht, d. h. ein Bild gegen ein anderes.

A4: Die Bilder sind bei Nachbestellungen teurer und im Display billiger als beim Kauf einzelner Päckchen.

Wir haben in Summe also drei der klassischen Annahmen angepasst (K1, K3, K4) und nur eine unverändert übernommen (K2). K5 ist zwar erfüllt, aber nicht praxisrelevant, sodass wir auf sie verzichten. Allerdings bedeuten unsere neuen praxisnahen Annahmen, dass wir dafür jetzt auch eine neue Sammelalbumformel finden müssen.

## 3   Die nützliche Sammelalbumformel

Wir führen zunächst die benötigten neuen Parameter ein, und zwar sei d die Anzahl der neuen Bilder in einem Display mit D Bildern und T die Anzahl der Karten, die der Sammler tauschen kann. Wir beginnen zunächst mit der Betrachtung des Nachkaufens.

**Satz 1:** Unter den Annahmen A1-A4 gilt für den Faktor

$$f = H(B - d) - H(K) + \frac{D+K}{B} \quad (2)$$

**Beweis:** Jedes Warten auf ein neues Bild bei x fehlenden Bildern kann als ein wiederholtes Zufallsexperiment mit Erfolgswahrscheinlichkeit $q = \frac{x}{B}$, dem Erwartungswert $\frac{1}{q}$ und der Varianz $\frac{1-q}{q^2}$ modelliert werden [6]. Der Erwartungswert der insgesamt benötigten Bilder ist die Summe der Erwartungswerte der einzelnen Wartezeiten. Auch die Varianzen addieren sich. Das Sammeln eines Albums kann also als eine Summe von unabhängigen, geometrisch verteilten Zufallsvariablen $X_1, X_2 \ldots X_B$ aufgefasst werden. Erhält man d Bilder beim Kauf eines Displays, so verringert sich der Faktor auf $H(B - d)$, da sich die Wartezeiten der letzten d Zufallsvariablen auf 1 verringern. Ebenso verringern sich durch das Nachkaufen die Wartezeiten der ersten K Zufallsvariablen ebenfalls auf 1 und der Faktor verringert sich um $H(K)$. Berücksichtigt man jetzt noch die zusätzlichen Doppelten und die Nachkaufbilder, die man einmal kaufen muss, so ergibt sich Formel (2).

Den Parameter d kann man ermitteln, indem man ein Display kauft und auspackt, oder indem man die Rezensionen bei Amazon auswertet [2].

Im nächsten Schritt erweitern wir unser Modell auf das Tauschen. Wir nehmen an, dass der Sammler eine feste Anzahl T von Bildern tauschen kann. Dabei muss T natürlich kleiner als B-d-K sein. Zur Vereinfachung machen wir eine wichtige Annahme, die jeder Sammler aus der Praxis bestätigen wird, nämlich

A5: Der Sammler hat immer genügend Bilder zum Tauschen.

Theoretisch könnte dies z. B. der Fall sein, wenn im Display nur neue Bilder wären, d. h. D=d oder wenn er beim Sammeln und Tauschen viel Glück hat. Normalerweise hat er am Anfang D-d Doppelte und wenn er nicht extrem konsequent tauscht, z. B. über Internet-Sammelbörsen, wird A5 in der Praxis immer erfüllt sein. A5 vereinfacht unser Model enorm, da wir sonst die aktuelle Anzahl der Doppelten berücksichtigen müssten.

Zunächst geben wir zwei einfache Abschätzungen

**Satz 2:** Unter den Annahmen A1-A5 gilt für den Faktor

$$H(B-d) - H(K+T) + \frac{D+K}{B} \leq f \leq H(B-d-T) - H(K) + \frac{D+K}{B} \quad (3)$$

**Beweis:** Die Summanden in der harmonischen Zahl $H(B)$ sind monoton fallend. Der Fall, dass man in (2) die ersten T Bilder gleich tauscht, entspricht einer Erhöhung der neuen Bilder im Display auf d+T und ergibt die obere Schranke. Der Fall, dass man in (2) die letzten T Bilder tauschen kann, entspricht einer Erhöhung der Nachkaufbilder auf K+T und ergibt die untere Schranke.

Unter der Bedingung, dass man genau T Bilder tauschen kann, entspricht dies der schlechtesten bzw. besten Tauschstrategie. Natürlich weiß man nicht, wann genau man die Bilder tauschen kann. Es ist daher relativ naheliegend, anzunehmen, dass das Tauschen quasi zufällig erfolgt. Dies bedeutet dass die Verteilung der Tauschbilder unter allen zu sammelnden Bildern quasi der Gleichverteilung entspricht.

**Satz 3:** Unter den Annahmen A1-A5 sowie der Annahme, dass T Bilder zu zufälligen Zeitpunkten getauscht werden, gilt für den Faktor

$$f = \ln\left(\left(\frac{B-d}{K}\right)^{1-t}\right) + \frac{D+K}{B} \quad (4)$$

wobei $t = \frac{T}{B-d-K}$.

**Beweis:** Die Annahme des zufälligen Tauschens bedeutet, dass in dem verbleibenden Teil von $H(B-d) - H(K)$ zufällig Summanden entfernt werden. Im Mittel bedeutet dies eine gleichmäßige Ausdünnung des Faktors zu $(H(B-d) - H(K))(1-t)$. Unter A5 sind immer genügend Tauschbilder vorhanden, d. h. es müssen keine Bilder extra gekauft werden, da alle Bilder gesammelt werden. Dann ergibt sich (4) einfach durch die Approximation von $H(B)$ durch den natürlichen Logarithmus.

Ein ähnliches Ergebnis können wir mit derselben Argumentation für die Standardabweichung erzielen.

**Satz 4:** Unter den Annahmen A1-A5 sowie der Annahme, dass T Bilder zu zufälligen Zeitpunkten getauscht werden, gilt

$$\frac{\sigma}{\sqrt{B}} = \sqrt{\left(\left(\frac{B-d}{K} - 1\right) - \ln\left(\frac{B-d}{K}\right)\right)(1-t)} \quad (5)$$

**Beweis:** Durch einfaches Einsetzen erhalten wir für die Varianz für einen Sammler ohne Nachkaufen

$$V = \sum_{i=1}^{B} \frac{B - \frac{i}{B}}{\frac{i^2}{B^2}} = B^2 \sum_{i=1}^{B} \frac{1}{i^2} - B \sum_{i=1}^{B} \frac{1}{i} = B^2 H_2(B) - B H(B)$$

wobei $H_2$ die verallgemeinerte harmonische Zahl zweiter Ordnung ist. Wir versuchen daher, für die Varianz des Einzelsammlers mit Nachkaufen dasselbe Argument wie bei dem Erwartungswert anzuwenden, nämlich dass die letzten K Terme in der Summe der Varianzen der Wartezeiten wegfallen. Damit erhalten wir

$$V = B^2(H_2(B) - H_2(K)) - B(H(B) - H(K)).$$

Jetzt können wir genauso die Anzahl der neuen Karten in Displays berücksichtigen und erhalten

$$V = B^2(H_2(B-d) - H_2(K)) - B(H(B-d) - H(K)).$$

Man erkennt jetzt, dass die Summanden in der Varianz ebenfalls monoton fallend sind und kann dasselbe Argument wie oben wiederholen. Unter Berücksichtigung der Approximation

$$H_2(n) \approx \int_1^n \frac{dx}{x^2} = 1 - \frac{1}{n}$$

erhalten wir schließlich (5).

Jetzt können wir die Ergebnisse für eine Formel für die Berechnung der mittleren Kosten des Sammelalbums zusammenfügen:

$$A = Bf\frac{p}{P} + (1+N)K\left(b - \frac{p}{P}\right) - \left(D\frac{p}{P} - C\right) \quad (6)$$

Dabei bezeichnet N die Anzahl der Nichtsammler, die zusätzlich für den Sammler mit nachkaufen sowie C die Kosten eines Displays. Der erste Summand entspricht den Kosten der Bilder, wenn sie alle gleich teuer wären, wobei f nach Formel 4 bestimmt wird. Die beiden nächsten Summanden berücksichtigen dann jeweils die Mehrkosten durch die Nachkaufkarten sowie die Kostenreduktion durch den Kauf des Displays.

## 4  Beispielberechnung für das EM-Album

Damit können wir die mittleren Kosten sowie deren Standardabweichung berechnen. Die Formeln haben wir in einem Open Office Tabellenkalkulationsblatt programmiert. Bei der Auswertung bemerken wir aber sofort, dass die Ergebnisse für große Tauschquoten t zu optimistisch sind, denn die Ergebnisse sind etwas unter dem Minimalpreis. Dies liegt daran, dass in diesen Fällen getauscht werden kann, ohne dass überhaupt genügend Doppelte vorhanden sind. Deshalb modifizieren wir die Ergebnisse so, dass immer mindestens der Minimalpreis gezahlt wird.

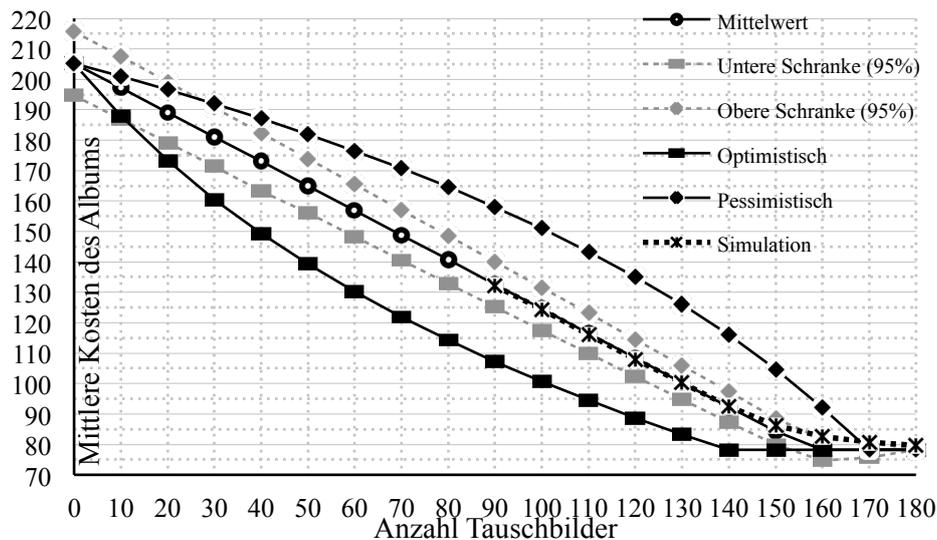

Abbildung 1: Beispielkosten für das EM-Album

Die Ergebnisse haben wir in Abbildung 1 für das EM-Album mit den Parametern d=450 und K=50 dargestellt, wobei wir zusätzlich noch die Kosten aufgetragen haben, die sich mit den pessimistischen bzw. optimistischen Abschätzungen nach Formel 3 ergeben würden. Außerdem haben wir den zweifachen Streuungsbereich mit dargestellt, d. h. den Bereich, in dem etwa 95% aller Sammler landen sollten.

In allen Fällen sieht man, dass die Ergebnisse für sehr hohe Tauschquoten t zu optimistisch sind, da sich ohne die Korrektur Gesamtkosten ergeben, die unter den minimal möglichen Kosten (hier 78,20 €) liegen würden. Dies liegt an unserer vereinfachenden Annahme, dass immer genügend Tauschbilder vorhanden sind. Diese Annahme ist im Fall sehr hoher Tauschquoten t nicht erfüllt und mit der Korrektur erreichen wir, dass wir mindestens beim Minimalpreis landen.

In Wirklichkeit muss man doch ein paar Bilder zusätzlich kaufen, sodass wir bei hohen Tauschquoten etwas zu optimistisch sind. Andererseits sieht man auch, wie stark es davon abhängt, wann man die Bilder tauscht und auch wie hoch die Tauschquote ist. Beides weiß man vorher noch nicht genau, d. h. die Unsicherheit über diese beiden Annahmen ist wahrscheinlich für die Prognose höher als der Fehler, den wir mit der vereinfachenden Annahme hereinbringen. Daher haben wir noch eine einfache Simulation in R geschrieben, bei welcher der Sammler nur tauschen kann, wenn er genug Doppelte hat. Die Kosten haben wir genauso mit der neuen Formel 6 berechnet. Die Simulation stimmt mit unserem berechneten Mittelwert gut überein. Erst ab 140 Tauschbildern weicht sie vom Mittelwert ab und verläuft etwa wie die obere Schranke (siehe Abbildung 1).

Schließlich wollen wir uns noch einem Phänomen widmen, das kurz nach dem Erscheinen des EM-Albums auftrat und allen unseren Berechnungen widersprochen hat.

## 5     Das Preis-Paradoxon

Schon kurz nach dem Erscheinen des EM-Sammelalbums gab es im Internet Angebote, deren Preis weit unter den von uns berechneten minimalen Kosten für die Sammelbilder lag. Die Angebote für das komplette Album lagen bei 120 € oder günstiger. Das könnte man nach unseren Untersuchungen nur mit sehr viel Tauschen und Nachkaufen schaffen. Tauschen wäre sehr aufwendig für jemanden, der damit viel Geld verdienen will. Für das Nachkaufen bräuchte man sehr viele Nichtsammler, was also auch eher unwahrscheinlich, da aufwändig, wäre. Manche Anbieter haben sogar Dutzende von Alben angeboten. Wenn wir annehmen, dass der Anbieter damit auch noch Gewinn machen möchte, dann müsste er das deutlich unter 100 Euro schaffen. Das ist paradox!

Wir schließen aus, dass der Anbieter illegal, also durch Diebstahl oder Firmenangehörige, an die Bilder kommt.

Zuerst wollen wir noch einmal die bisherigen Fakten, die uns bekannt sind, zusammenstellen:

1. Wir haben nachgewiesen [2], dass Panini die Bilder schlecht mischt.

2. Wenn die Bilder nacheinander aus derselben Fifimatic kämen, dann würden nacheinander die Bilder eines kompletten Albums kommen [2].

3. Wir haben bei Twitter in einer Meldung [13] gesehen, dass im Großhandel immer zwölf Displays zu einer „Stange" (so ähnlich wie bei Zigaretten) verpackt werden. Diese Stangen werden dann zu größeren Gebinden zusammengepackt, die auf einer Euro-Palette ausgeliefert werden.

Wir wissen nicht, wie die Displays auf die Stangen verteilt werden. Daher unterscheiden wir zwei mögliche Fälle:

Fall 1: In einer Stange sind nur Displays aus einer Fifimatic.

12 Displays entsprechen 12x500 Bildern, d. h. 6000 Bildern. Damit könnte man 8,8 Alben füllen, wenn man keine Doppelten hätte. Da wir wissen, dass Panini schlecht mischt, nehmen wir an, dass der Verkäufer damit wirklich 8 Alben füllen kann.

Ein Display kostete bei Amazon etwa 50 €, eine Stange würde also 600 € kosten (im Großhandel wahrscheinlich noch weniger). Damit würde jedes Album weniger als 75 € kosten und das Paradoxon kann damit erklärt werden. Gegenüber einem Verkaufspreis von 120 € könnte man ein gutes Geschäft machen. Der Verkäufer muss allerdings alle Päckchen auspacken und die Bilder sortieren.

Fall 2: Die Displays in einer Stange kommen nicht immer aus derselben Fifimatic.

Der Verkäufer müsste viele Stangen kaufen und auf die Seriennummern achten. Auf den Stangen sind mehrere Nummern angegeben [13], da aber nur eine Stange fotografiert wurde, wissen wir nicht, welche der Nummern die Seriennummer ist. Allerdings ist das Produktionsdatum sowie die Uhrzeit erkennbar. Damit könnte er hoffen, dass er wieder Bilder aus aufeinanderfolgenden Displays bekommt, die allerdings nicht nur aus einer Fifimatic stammen und er müsste viel mehr Stangen kaufen. Dies würde ebenfalls den Preis erklären, aber der Verkäufer würde ein viel höheres Risiko eingehen.

## 6    Diskussion und Zusammenfassung

Wir haben in dieser Arbeit neue Annahmen definiert und begründet, die zu realistischeren Ergebnissen und für den Sammler wesentlich günstigeren Prognosen, als beim klassischen Modell führen. Wir haben bisher nur eine Arbeit gefunden [12], in der von den klassischen Annahmen abgewichen wurde. Dabei haben wir Annahmen getroffen, die wie A1 neu sind oder wie A3 bisher zumindest selten getroffen wurden, aber die heutige Sammelpraxis gut widerspiegeln. Natürlich muss man einräumen, dass sich seit den 1950er Jahren, als die klassischen Annahmen getroffen wurden [5], vieles geändert hat, so z. B. auch das Tauschverhalten, beispielsweise über das Internet.

Von praktischer Bedeutung sind insbesondere unsere Formeln 4 und 5, mit denen man die mittleren Kosten des komplettierten Albums sowie die Standardabweichung der Kosten unter den neuen Annahmen ausrechnen kann. Wir meinen, dass wir damit die erste realistische, allgemeingültige Lösung für das Sammelbilderproblem gefunden haben. Sie ist nach unserer Literaturkenntnis, aber auch nach Wikipedia [1], für diese Kombination der Annahmen neu. Dies wird ebenfalls dadurch gestützt, dass die Ergebnisse, die kürzlich in den Medien diskutiert wurden [7, 9], auf den klassischen Annahmen aufbauen.

Zwar wurden schon Ergebnisse unter ähnlichen Annahmen wie A3 erzielt [12], aber die Annahmen zum Tauschen wurden nicht explizit ausgewiesen. Dort wurde die Wahrscheinlichkeitsverteilung mit einem Baumdiagramm berechnet oder simuliert, aber es wurden keine Formeln abgeleitet. Dies ähnelt den Ansätzen in unserer ersten Arbeit [5]. Die Verteilung ist allerdings auch schon vorher bekannt gewesen [6]. Aus den Ergebnissen kann man schließen, dass es sich für das Tauschen wohl eher um optimistische Annahmen handelt wie bei unserer Formel 14. Es wird nämlich geraten, „solange Päckchen zu kaufen, bis man … entsprechendes Tauschmaterial zusammen hat. Danach tauscht man … und bestellt die dann noch fehlenden Sticker nach."

In Tabelle 1 haben wir die Formeln für den Faktor noch einmal übersichtlich gegenübergestellt. Für das EM-Album ist der Faktor für das klassische Ergebnis nach Formel 2 etwa 7,1. Bei einer sehr großen Sammelgemeinschaft sinkt der Faktor nach Formel 5 auf ungefähr 1,88, aber nicht auf 1. In unserer letzten Arbeit [4] hatten wir schon gezeigt, dass der Effekt des Nachkaufens groß ist. Für den Einzelsammler sinkt der Faktor alleine dadurch auf etwa 2,68. Per Simulation konnten wir zeigen, dass der Faktor bei einer großen Sammelgemeinschaft gegen 1 geht. Unsere neuen Ergebnisse zeigen, dass dies unter unserer neuen Annahme A3 schon für den Einzelsammler ohne Nachkaufen gilt, wenn die Tauschquote gegen 1 geht. Nehmen wir mal an, dass der Einzelsammler nur 100 Bilder tauschen kann, dann sinkt der Faktor schon auf 5,56. Wenn er dann noch nachkauft, sinkt der Faktor auf 2,3. Im Vergleich dazu ist der Effekt des Display-Kaufs eher gering, denn selbst wenn es im Display gar keine Doppelten gäbe, wäre der Faktor immer noch 5,92. Nimmt man an, dass im Display 100 Doppelte wären, so wäre der Faktor sogar 6,37.

| Sammelbilderproblem | Ohne Nachkaufen | Mit Nachkaufen | und Display |
|---|---|---|---|
| **Einzelsammler** | $\ln B$ | $\ln \frac{B}{K} + \frac{K}{B}$ | $\ln \frac{B-d}{K} + \frac{K+D}{B}$ |
| **Tauschen (K3)** für großes F | $\ln \ln B$ | Gegen 1 | Gegen 1 |
| **Tauschen (A3)** | $\ln(B^{1-t})$ | $\ln\left(\left(\frac{B}{K}\right)^{1-t}\right) + \frac{K}{B}$ | $\ln\left(\left(\frac{B-d}{K}\right)^{1-t}\right) + \frac{D+K}{B}$ |

**Tab. 1.** Asymptotik des Faktors bei verschiedenen Varianten des Sammelbilderproblems

Wir müssen allerdings einräumen, dass unsere Formeln auch mit der Korrektur für sehr hohe Tauschquoten nicht ganz korrekt sind, d. h. wenn t gegen 1 geht. Dies liegt an der vereinfachenden Annahme, dass immer genügend Tauschbilder vorhanden sein sollen. Sonst müssten wir irgendwie noch Buch führen, wie viele Doppelte der Sammler wirklich hat, wann er tauscht usw. Dazu bräuchte man aber ein viel komplizierteres Modell, z. B. für die anderen Tauschpartner. Wir haben aber per Simulation für ein typisches Beispiel gezeigt, dass die Abweichungen gering sind und das Ergebnis stärker von der Tauschstrategie abhängt als von unserer vereinfachenden Annahme, sodass wir die Genauigkeit der Ergebnisse zumindest für praktische Zwecke für ausreichend erachten.

Außerdem könnte man kritisieren, dass wir angenommen haben, dass die Zahl der neuen Bilder im Display d konstant ist. Das bedeutet, dass man die Rechnung macht, nachdem man das Display ausgepackt hat. Will man vorher eine Prognose machen, müsste man noch die Standardabweichung von d berücksichtigen. Wir haben die Amazon-Rezensionen ausgewertet [10]. Dort war die Standardabweichung mit ungefähr 30 zwar groß, aber der Unterschied, den das ausmacht, ist relativ gering, da die einzelnen Beiträge zur Standardabweichung klein sind. Alternativ könnte man für d eine Verteilung annehmen, und man müsste dann die Standardabweichung ausrechnen oder simulieren.

Schließlich bleibt als größtes Problem für den Anwender unserer Ergebnisse die Einschätzung der Anzahl T der Tauschbilder bzw. der Tauschquote t. Dies mag für einen Anfänger schwierig sein, ein erfahrener Sammler sollte allerdings schon eine grobe Abschätzung machen können. Auch hier gilt wieder, dass sich geringe Fehleinschätzungen nicht stark auswirken, und man im Zweifel auch noch unsere optimistische bzw. pessimistische Abschätzung zu Rate ziehen kann (bzw. seine eigene Tauschstra-tegie dahingehend optimieren kann), denn wir berechnen die mittleren Kosten unter der Annahme, dass zufällig, d. h. spontan getauscht wird.

Wir meinen, dass wir ein für die Praxis sehr nützliches Modell gefunden haben und dass der Aufwand, es noch weiter zu verfeinern, sehr hoch bzw. aus unserer Sicht unangemessen wäre. Für uns stellt dieses Ergebnis daher den Abschluss unserer Forschung zu dem Thema dar, da unserer Meinung nach alle relevanten Forschungsfragestellungen befriedigend beantwortet sind. Außerdem sind alle Modelle ja nur Annäherungen an die Realität, wie es auch das berühmte Zitat „Essentially, all models are wrong, but some are useful" von George Box treffend ausdrückt.

## Referenzen

1. Wikipedia: *Sammelbilderproblem*, http://de.wikipedia.org/wiki/Sammelbilderproblem, letzter Abruf 17.2.2017
2. Braband, N., Braband, S., Braband, M.: *Das Geheimnis der Fifimatic- oder: Betrügen Sammelbildhersteller?*, erscheint in: Junge Wissenschaft, Ausgabe Nr. 114, 2017 (Preprint: arXiv:1603.03008, 2016)
3. Sardy, S., Velenik, Y.: *Paninimania: Sticker rarity and cost-effective strategy*, Université de Genève (2014), http://www.unige.ch/math/folks/velenik/Vulg/Paninimania.pdf
4. Braband, N., Braband, S.: *Nicht mehr über Sammelbilder ärgern!*, Junge Wissenschaft, Ausgabe Nr. 110, 2016, S.16-24
5. Newman, L., Shepp, L.: *The double dixie cup problem*, The American Mathematical Monthly (1960), 58-61, http://www.jstor.org/stable/2308930
6. Henze, N.: *Stochastik für Einsteiger*. 11. Auflage. Springer Spektrum, Wiesbaden 2016.
7. Wales Online: *A maths genius worked out exactly how much it will cost to fill your Panini Euro 2016 album.* http://www.walesonline.co.uk/sport/football/football-news/maths-genius-worked-out-exactly-11120318 (letzter Abruf 4.5.2016).
8. Private email-Kommunikation mit Prof. Paul Harper (Cardiff University), 27.4.-3.5.2016
9. Deutsche Mathematiker Vereinigung: *Panini-Bilder, Mathe-Genies und Jugend forscht*, DMV Forum, https://dmv.mathematik.de/index.php/forum/nachrichten/533-panini-bilder-mathe-genies-und-jugend-forscht (letzter Abruf 11.12.2016)
10. Amazon: *Panini Euro 2016 Stickerbox mit 100 Tüten X 5 Sticker = 500 Sticker*, ASIN B01CPN3NYO (letzter Abruf 17.12.2016)
11. DER SPIEGEL: *Da fehlen doch welche*, Heft 25/2014, S. 78, auch: http://www.spiegel.de/spiegel/print/d-127626366.html (letzter Abruf 1.2.2016)
12. Binzenhöfer, A., Hoßfeld, T.: *Warum Panini Fußballalben auch Informatikern Spaß machen*. In: Hans-Georg Weigand (Hrsg.): Fußball – eine Wissenschaft für sich. Königshausen & Neumann, Würzburg 2006, S. 181–191.
13. Noel Clarke: *Found the countries supply of panini #EURO2016 stickers…*, https://twitter.com/NoelClarke/status/728710989609324544 (letzter Abruf 11.12.2016)